\documentclass[a4paper]{article}
\usepackage{amssymb,amsmath,amsthm,verbatim,IEEEtrantools,macros,mathtools}
\usepackage[margin=1in]{geometry}
\DeclareMathOperator{\Exact}{Exact}
\DeclareMathOperator{\Tight}{Tight}
\DeclareMathOperator{\supp}{supp}
\DeclareMathOperator{\dimf}{dim_F}
\title{Fourier Dimension Estimates for Sets of Exact Approximation Order: the Well-Approximable Case}
\author{Robert Fraser \and Reuben Wheeler}
\begin{document}
\maketitle
\begin{abstract}
We obtain a Fourier dimension estimate for sets of exact approximation order introduced by Bugeaud for certain approximation functions $\psi$. This Fourier dimension estimate implies that these sets of exact approximation order contain normal numbers.
\end{abstract}
\section{Introduction and Background}
\subsection{Hausdorff and Fourier Dimension}
Let $E \subset \mathbb{R}$ be a compact set. Frostman's lemma \cite{Frostman35} implies that the Hausdorff dimension $\text{dim}_H(E)$ of $E$ is the supremum over all values of $s < 1$ such that $E$ supports a Borel probability measure $\mu$ satisfying the condition
\[\int |\hat \mu(\xi)|^2 |\xi|^{s-1} d \xi < \infty.\]
This condition essentially says that $|\hat \mu(\xi)|$ decays, in an $L^2$-average sense, like $|\xi|^{-s/2}$. 

If, instead, we impose the condition that $|\hat \mu(\xi)|$ decays \textit{pointwise} like $|\xi|^{-s/2}$ for all $s < s_0 \leq 1$, we say that $E$ has Fourier dimension at least $s_0$, written $\text{dim}_F(E) \geq s_0$. Clearly, we have that $\text{dim}_F(E) \leq \text{dim}_H(E)$ for every compact subset $E \subset \mathbb{R}$. If a compact set $E$ satisfies $\text{dim}_F(E) = \text{dim}_H(E)$, then we say that $E$ is a compact \textbf{Salem set}.

In fact, it is nontrivial to construct examples of Salem sets. The earliest constructions of Salem sets are random constructions, such as random cantor sets of Salem \cite{Salem51}. K\"orner established via a random construction that the Fourier dimension of a set $E \subset \mathbb{R}$ can take any value from $0$ to $\text{dim}_H(E)$ \cite{Korner11}.

\subsection{Metric Diophantine Approximation}
A classical result of Jarn\`ik and Besicovitch \cite{Besicovitch34} \cite{Jarnik29} concerns the Hausdorff dimension of the set of $\tau$-approximable numbers. The $\tau$-approximable numbers are the set
\[E(\tau) := \{x : |x - r/q| \leq q^{-\tau} \, \text{for infinitely many pairs of integers $(r,q)$.}\}\]
For $\tau \leq 2$, it is easy to see using the Dirichlet principle that $E(\tau) = \mathbb{R}$. Jarn\`ik and Besicovitch show that for $\tau > 2$, $\text{dim}_H(E(\tau)) = \frac{2}{\tau}$.

Kaufman \cite{Kaufman81} established that, in fact, the set $E(\tau)$ has Fourier dimension equal to $\frac{2}{\tau}$, implying that $E(\tau)$ is a Salem set. Notably, this is the first explicit non-random construction of a Salem set of Hausdorff dimension other than $0$ or $1$ in $\mathbb{R}$.

Fourier dimension calculations are of interest in metric Diophantine approximation because of a celebrated result of Davenport, Erd\"os, and Leveque \cite{DavenportErdosLeveque63}. This result concerns the presence of normal numbers in subsets $E \subset \mathbb{R}$. 

Specifically, Davenport, Erd\"os and Leveque show that if $\mu$ is a positive Borel probability measure, and $a$ is a positive integer, then the sequence $\{a^jx\}_{n=1}^{\infty}$ is uniformly distributed modulo $1$ if and only if
\begin{equation}\label{DELequation1}
\sum_{N = 1}^{\infty} N^{-3} \sum_{j=1}^N \sum_{k=1}^N \hat \mu(m(a^j - a^k)) < \infty
\end{equation}
for every nonzero integer $m \in \mathbb{Z}$.  If we crudely assume $|\hat \mu(\xi)| \leq |\xi|^{-s/2}$ for every $\xi \in \mathbb{R}$, the sum in $k$ in \eqref{DELequation1} is essentially a geometric sum, so we have that the sum in $j$ is bounded above by an $m$-dependent constant times $N$, and the sum certainly converges for all nonzero integers $m$, and all integers $a \geq 2$. Therefore, if $\mu$ is a Borel probability measure such that $|\hat \mu(\xi)| \leq |\xi|^{-s/2}$ for some $s > 0$, then $\mu$-almost every point is a normal number. Note that a far weaker assumption on $\mu$ suffices to locate normal numbers; see e.g. \cite{PollingtonVelaniZafeiropoulosZorin19}.

In particular, any set $E$ of positive Fourier dimension must contain normal numbers. It is therefore of interest to find Fourier dimension estimates for subsets of $\mathbb{R}$ arising in Diophantine approximation. Of course, the well-approximable numbers, being a Salem set of positive dimension, contain normal numbers. In fact, Kaufman also showed a Fourier dimension result for sets of \textit{badly}-approximable numbers \cite{Kaufman80}.

The badly approximable numbers consist of those $x \in \mathbb{R}$ such that the partial quotients in the continued fraction expansion of $x$ are bounded. Given a finite set $S \subset \mathbb{N}$ with at least two elements, we use the term $S$-badly-approximable numbers to refer to those real numbers $x$ such that the partial quotients of the continued fraction expansion of $x$ all lie in the finite set $S$. Kaufman \cite{Kaufman80} established that, if $S$ is a finite set such that the Hausdorff dimension of the $S$-badly-approximable numbers is greater than $2/3$, then the $S$-badly approximable numbers have positive Fourier dimension.

The method used by Kaufman to estimate the Fourier dimension of the badly approximable numbers is very different from the method used to estimate the Fourier dimension of the well-approximable numbers. For the well-approximable numbers, Kaufman's argument relies on the cancellation of the exponential sum
\begin{equation}\label{geomseries}
\sum_{r=0}^{q} e(rs/q)
\end{equation}
for any integers $s, q$ such that $q$ does not divide $s$; since $s$ has a small number of divisors, the sum of \eqref{geomseries} over all $M/2 \leq q < M$ will also be small.

In contrast, Kaufman's Fourier dimension estimate \cite{Kaufman80} follows a rather different argument. This argument relies on constructing a certain random measure on bounded integer sequences whose pushforward under the continued fraction map satisfies the relevant Fourier decay condition, which is established via a van der Corput-type lemma.

Queffelec and Ramar\'e \cite{QueffelecRamare03} improve the $2/3$ requirement on the Hausdorff dimension to $1/2$; in particular, this condition holds if $S = \{1,2\}$. Hochman and Shmerkin \cite{HochmanShmerkin15} show that, without any Hausdorff dimension assumption, the set of $S$-badly-approximable numbers contains normal numbers for any finite set $S \subset \mathbb{N}$ with at least two elements. In a recent work, Sahlsten and Stevens \cite{SahlstenStevens20} improved on all of these results by showing, without any Hausdorff dimension assumption, that the $S$-badly-approximable numbers have positive Fourier dimension for any finite set $S \subset \mathbb{N}$ with at least two elements.

Of note is that, while Kaufman's argument for the well-approximable numbers \cite{Kaufman81} works just as well in the inhomogeneous setting, there does not seem to be an easy way to modify Kaufman's argument for the badly approximable numbers, \cite{Kaufman80}, to this case. Doing so would require a satisfactory analogue of the continued fraction expansion for the inhomogeneous version of the badly approximable numbers.
\subsection{Approximation to Exact Order}
Bugeaud \cite{Bugeaud03} introduced sets of exact approximation order. We will now define an inhomogeneous analogue. Given an approximation function $\psi$ and a real number $\theta \in [0,1)$, define the set $\Exact(\psi, \theta)$ to be the set of real numbers $x$ satisfying the pair of conditions:
\begin{IEEEeqnarray*}{rCLl}
\left|x - \frac{r - \theta}{q} \right| & \leq & \psi(q) \quad & \text{for infinitely many pairs $(r,q)$ of relatively prime integers} \\
\left|x - \frac{r - \theta}{q} \right| & \leq & \psi(q) - c \psi(q) \quad & \text{for only finitely many pairs $(r,q)$ of relatively prime integers and any $c > 0 $.} \\
\end{IEEEeqnarray*}
Bugeaud \cite{Bugeaud03} computes the Hausdorff dimension of the set $\Exact(\psi, 0)$ for certain functions $\psi$. Specifically, Bugeaud considers functions $\psi$ such that the function $x^2 \psi(x)$ is nonincreasing. Bugeaud shows that the Hausdorff dimension of $\Exact(\psi, 0)$ is $\frac{2}{\lambda^+}$, where 
\[\lambda^+(\psi) = - \limsup_{x \to \infty} \frac{\log \psi(q)}{\log q} \geq 2.\]
The upper Hausdorff dimension bound follows trivially from the Jarn\`ik-Besicovitch theorem. For the lower bound, Bugeaud considers a subset of $\Exact(\psi,0)$ consisting of numbers whose continued fractions have partial quotients that typically grow very slowly, except for some exceptional partial quotients that are larger.

The main result in this paper concerns the Fourier dimension of the sets $\Exact(\psi, \theta)$. Unlike Bugeaud, we will limit ourselves to approximation functions $\psi$ of satisfying the property that
\begin{equation}\label{orderdef}
\lambda(\psi): =- \lim_{q \to \infty} \frac{\log \psi(q)}{\log q}
\end{equation}
exists. 

Observe that the set $\Exact(\psi, \theta)$ is invariant under translations by integers. Therefore, we can view $\Exact(\psi, \theta)$ as a subset of the torus $[0,1)$ in a natural way. We will use the notation $\Exact^{[0,1)}(\psi, \theta)$ to refer to this subset of the torus.

Let $\theta \in [0,1)$ be an irrational number. We define the \textbf{Diophantine approximation exponent} of $\theta$ to be the infimum over values of $\gamma$ such that the equation
\[\left|\theta - \frac{p}{q} \right| \leq q^{-\gamma}\]
has only finitely many solutions for integers $p$ and $q$. Note that this implies that
\[\left|q \theta - p \right| \leq q^{-\gamma + 1}\]
also has only finitely many solutions with $q \neq 0$.

This means that if $\gamma$ is the Diophantine approximation exponent of $\theta$, $\eta > 0$, then we have that for any integers $p$ and $q$,
\[\left|q \theta - p \right| \geq q^{-\gamma - \eta + 1},\]
provided $q$ is sufficiently large depending on $\eta$ and $\theta$. 
Another way of saying this is that
\[\norm{q \theta} \geq q^{- \gamma - \eta + 1}\]
where $\norm{q \theta}$ is the distance from $q \theta$ to the nearest integer.
\begin{mythm}\label{mainthm}
Let $\theta \in [0,1)$ be either $0$ or an irrational number with finite Diophantine exponent $\gamma$, taking $\gamma = 1$ if $\theta = 0$. Let $\psi(q)$ be a positive, decreasing function of $q$. Suppose $\tau = \lambda(\psi)$ is such that
\begin{equation}\label{taucondition}
\tau > \frac{2 + \gamma^2 + \sqrt{(2 + \gamma^2)^2 - 4}}{2}.
\end{equation}
Then $\dimf \Exact(\psi, \theta)$ is positive; moreover, we have the inequality
\begin{equation}\label{fourierdimestimate}
\dimf \Exact(\psi, \theta) \geq \alpha := \frac{2(\beta - \tau)}{\tau (\beta - 1)},
\end{equation}
where
\begin{equation}\label{betadef}
\beta := \gamma^{-2} (\tau - 1)^2.
\end{equation}
\end{mythm}
One can observe that the condition \eqref{taucondition} implies that $\beta > \tau$; this implies that the right side of \eqref{fourierdimestimate} is positive. We will make some quick observations about Theorem \ref{mainthm}.
\begin{myrmk}
If $\theta = 0$, we are able to take $\gamma = 1$. In this case, the inequality \eqref{taucondition} reduces to
\begin{equation}\label{tauconditionhomog}
\tau > \frac{3 + \sqrt{5}}{2}.
\end{equation}
\end{myrmk}
\begin{myrmk}
For a fixed $\gamma$, observe that in the regime $\tau \to \infty$, we have that $\beta = \gamma^{-2} \tau^2 + O(\tau)$. Therefore,
\[\lim_{\tau \to \infty} \dfrac{\frac{2(\beta - \tau)}{\tau (\beta - 1)}}{2/\tau} = 1.\]
This means that, if $\tau$ is large, Theorem \ref{mainthm} ``nearly" shows the set $\Exact(\psi, \theta)$ is a Salem set.
\end{myrmk}

In fact, our proof yields a slightly more general Fourier dimension estimate than the one in Theorem \ref{mainthm}. In order to state this estimate, we will introduce sets of tight approximation order.

\begin{mydef}[Sets of tight approximation order]
Let $\psi_1$, $\psi_2$ be a pair of approximation functions with $\psi_2(q) \leq \psi_1(q)$ for all $q$. The set $\Tight(\psi_1, \psi_2, \theta)$ consists of those real numbers $x$ satisfying the conditions
\begin{IEEEeqnarray*}{rCLl}
\left|x - \frac{r - \theta}{q} \right| & \leq & \psi_1(q) \quad & \text{for infinitely many pairs $(r,q)$ of relatively prime integers,} \\
\left|x - \frac{r - \theta}{q} \right| & \leq & \psi_1(q) - c \psi_2(q) \quad & \text{for only finitely many $(r,q)$ of relatively prime integers and any $c > 0$.} \\
\end{IEEEeqnarray*}
\end{mydef}
The set $\Exact(\psi, \theta)$ is the same as the set $\Tight(\psi, \psi, \theta)$. Given appropriate conditions on $\psi_1$ and $\psi_2$, we are able to estimate the Fourier dimension of the set $\Tight(\psi_1, \psi_2, \theta)$.
\begin{mythm}\label{mainthm2}
Let $\theta \in [0,1)$ be either $0$ or an irrational number with finite Diophantine exponent $\gamma$, taking $\gamma = 1$ if $\theta = 0$. Let $\psi_1(q)$ and $\psi_2(q)$ be decreasing functions with $\psi_2(q) \leq \psi_1(q)$ for all $q$ and such that $\psi_1(q) - \psi_2(q)$ is decreasing. 

Let

\begin{equation}\label{betadef2}
\beta := \gamma^{-2} (\tau_1 - 1)^2.
\end{equation}

Suppose $\tau_1 = \lambda(\psi_1)$ and $\tau_2 = \lambda(\psi_2)$ are such that $\beta > \tau_2$. 
Then $\dimf \Tight(\psi_1, \psi_2, \theta)$ is positive; moreover, we have the inequality
\begin{equation}\label{fourierdimestimate2}
\dimf \Tight(\psi_1, \psi_2, \theta) \geq \alpha := \frac{2(\beta - \tau_2)}{\tau_2 (\beta - 1)}.
\end{equation}
\end{mythm}
\section{An elementary Diophantine approximation lemma}
The key to adapting Kaufman's argument to the set $\Tight(\psi_1, \psi_2, \theta)$ is an elementary lemma in Diophantine approximation. This lemma states that if a real number $x$ is approximable by rationals at two ``fairly close" scales, then $x$ cannot be approximable at any intermediate scale.

Let \begin{equation}\label{eq:betaEps}\beta_{\epsilon} = \gamma^{-2}(\tau_1- 1)^2 - \epsilon.\end{equation}
\begin{mylem}\label{DiophantineLemmaGeneral}
Let $\tau_1 > 2$, $0 < c < 1$, $\epsilon > 0$ be real numbers, and let $x \in \mathbb{R}$. Suppose $\psi_1(q)$ is a decreasing function with the property that $\lambda(\psi_1) = \tau_1$ and $\psi_2(q)$ is chosen with $\psi_2$ and $\psi_1 - \psi_2$ decreasing, such that $\lambda(\psi_1) = \tau_1$, and $\lambda(\psi_2) = \tau_2$. Let $\theta \in [0,1)$ be either $0$ or an irrational number with Diophantine approximation exponent $\gamma$, taking $\gamma = 1$ if $\theta = 0$.  If $x$ satisfies the pair of inequalities:
\begin{IEEEeqnarray*}{rCCCl}
\psi_1(q_1) - c \psi_2(q_1) & \leq & \left|x - \frac{p_1 - \theta}{q_1} \right| & \leq & \psi_1(q_1)\\
\psi_1(q_2) - c \psi_2(q_2) & \leq & \left|x - \frac{p_2 - \theta}{q_2} \right| & \leq & \psi_1(q_2)
\end{IEEEeqnarray*}
where $Q(\epsilon) < q_1 <  q_2  < q_1^{\beta_{\epsilon}}$ and $q_2$ is prime, then $x$ does not satisfy any inequality of the form
\[\left|x - \frac{p - \theta}{q} \right| < \psi_1(q) - c \psi_2(q)\]
for $q_1 < q < q_2.$
\end{mylem}
\begin{proof}
Suppose $x$ satisfies the conditions of Lemma \ref{DiophantineLemmaGeneral}. Then we have the inequality
\begin{equation}\label{q1ineq}
\psi_1(q_1) - c \psi_2(q_1) \leq \left|x - \frac{p_1 - \theta}{q_1} \right| \leq \psi_1(q_1)\end{equation}
We will now split into two cases depending on whether $\theta = 0$.
\paragraph*{Case 1} Here we consider $\theta = 0$. In this case, \eqref{q1ineq} reduces to
\begin{equation}\label{q1ineqzero}
\psi_1(q_1) - c \psi_2(q_1) \leq \left|x - \frac{p_1}{q_1} \right| \leq \psi_1(q_1).
\end{equation}
Observe that if there exist $p, q$ such that 
\begin{equation}\label{falseassumption}
\left|x - \frac{p}{q} \right| < \psi_1(q) - c \psi_2(q),
\end{equation}
with $q > q_1$, then we must have $\frac{p}{q} \neq \frac{p_1}{q_1}$; this follows from \eqref{q1ineqzero} and the fact that $\psi_1 - c \psi_2$ is decreasing.

Thus, we have the inequality
\[\left| \frac{p}{q} - \frac{p_1}{q_1} \right| = \frac{|pq_1 - p_1q|}{q q_1} \geq \frac{1}{q q_1}\]
since the numerator is a nonzero integer. On the other hand, equations \eqref{q1ineqzero} and \eqref{falseassumption} imply by the triangle inequality that
\[\left| \frac{p}{q} - \frac{p_1}{q_1} \right| \leq \psi_1(q_1) + \psi_1(q) - c \psi_2(q) \leq 2 \psi_1(q_1),\]
where the last inequality follows from the fact that $\psi_1$ is decreasing.

Combining these inequalities gives that 
\begin{equation}\label{qq1ineq}
\frac{1}{q q_1} \leq 2 \psi_1(q_1).
\end{equation}
Take $\eta > 0$. At this stage, we select some $Q_{\eta}$ such that $\frac{\log (2\psi_1(q))}{\log q} < - \tau_1 + \eta$ for all $q \geq Q_\eta$. Then, it follows from \eqref{qq1ineq}, that if $q_1 \geq Q_{\eta}$,
\[q \geq q_1^{\tau_1 - 1 - \eta}.\]
A similar argument, using the fact that $\frac{p}{q} \neq \frac{p_2}{q_2}$ by the primality of $q_2$, reveals that 
\[q_2 \geq q^{\tau_1 - 1 - \eta}\]
and thus, we combine to get that
\[q_2 \geq  q_1^{(\tau_1 - 1 - \eta)^2}.\]
So if we choose $\eta$ sufficiently small that $(\tau_1 - 1 -\eta)^2 > (\tau_1 - 1)^2 - \epsilon$, we must have
\[q_2 \geq q_1^{\beta_{\epsilon}},\]
as desired.
\paragraph*{Case 2} Now, we will assume $\theta \in [0,1)$ is an irrational number with Diophantine approximation exponent $\gamma$. In this case, we simply observe that if \eqref{q1ineq} and an analogue of \eqref{falseassumption} both hold, then we must have, by the triangle inequality and the fact that $\psi_1$ is decreasing, that
\begin{IEEEeqnarray}{rCl}
\left|\frac{p_1 - \theta}{q_1} - \frac{p - \theta}{q} \right| & \leq & 2 \psi_1(q_1) \\
\text{and}\qquad\left|\frac{q p_1 - p q_1 - (q - q_1)\theta}{q q_1} \right| & \leq & 2 \psi_1(q_1).\label{q1ineqinhom}
\end{IEEEeqnarray}
Now, $q p_1 - p q_1$ is an integer, as is $(q - q_1)$. Therefore, if $\theta$ has Diophantine exponent $\gamma$, then we have $||(q - q_1) \theta|| \geq C_{\eta, \theta}(q - q_1)^{- \gamma + 1- \eta/2} \geq q^{-\gamma + 1 - \eta}$, provided that $q > Q_{\eta,1}$, where $Q_{\eta,1}$ is sufficiently large. Here, $||\cdot||$ denotes the distance to the nearest integer. Thus, for such $q$, we have
\begin{equation}\label{falseassumptioninhom}
\left |\frac{q p_1 - p q_1 - (q - q_1)\theta}{q q_1} \right| \geq \frac{q^{- \gamma + 1 - \eta}}{q q_1}.
\end{equation}
By combining inequalities \eqref{q1ineqinhom} and \eqref{falseassumptioninhom}, we get
\[q^{- \gamma - \eta} q_1^{-1} \leq 2 \psi_1(q_1).\]
Now, observe that if $q_1 \geq Q_{\eta, 2}$, then we have that $2 \psi_1(q_1) \leq q_1^{-\tau_1 + \eta}$ as in Case 1. Thus, for $q_1 > \max(Q_{\eta, 1}, Q_{\eta, 2})$, we then have
\[q^{-\gamma - \eta} q_1^{-1} \leq q_1^{- \tau_1 + \eta}\]
and solving for $q$ yields
\[q \geq q_1^{(\gamma + \eta)^{-1}(\tau_1 - 1 - \eta)}.\]
By a similar argument, we observe
\[q_2 \geq q^{(\gamma + \eta)^{-1}(\tau_1 - 1 - \eta)}.\]
Combining these inequalities gives
\[q_2 \geq q_1^{(\gamma + \eta)^{-2} (\tau_1 - 1 - \eta)^2}.\]
If $\eta$ is sufficiently small relative to $\epsilon$, we then have the inequality
\[q_2 \geq q_1^{\beta_{\epsilon}}\]
as desired.
\end{proof}

For the purposes of the rest of the argument, it will be important to have $\beta_{\epsilon} > \tau_2$ for sufficiently small $\epsilon$. This leads to the restriction
\eqref{taucondition}.
\section{A periodization trick}
In order to establish Theorem \ref{mainthm2}, we must construct, for any $\epsilon > 0$, a finite Borel measure $\mu_{\epsilon}$ supported on $\Tight(\psi_1, \psi_2, \theta) \cap [0,1]$ such that $\hat \mu_{\epsilon}(\xi) \leq C |\xi|^{-\alpha/2 + \epsilon}$ for all $\xi \in \mathbb{R}$. However, it is convenient to evaluate $\hat \mu_{\epsilon}(\xi)$ at only integer values of $\xi$. For this purpose, it is convenient to introduce the set $\Tight^{[0,1)}(\psi_1, \psi_2, \theta)$, a subset of the torus. To this end, we will construct a measure $\mu_{\epsilon}^{[0,1)}$ on the torus with support contained in $\Tight^{[0,1)}(\psi_1, \psi_2, \theta)$. Observe that, as $\mu_{\epsilon}^{[0,1)}$ is a measure on the torus, it has a corresponding Fourier-Stieltjes series, the coefficients of which will be denoted $\widehat{\mu_{\epsilon}^{[0,1)}}(\xi)$. Such a measure $\mu_{\epsilon}^{[0,1)}$ can be associated to a $1$-periodic measure $\mu_{\epsilon}^P$ supported on the real numbers.

\begin{mylem}\label{torusmeasure}
Suppose that $\mu_{\epsilon}^{[0,1)}$ is a measure on the torus with support contained in $\Tight^{[0,1)}(\psi_1, \psi_2, \theta)$, with the property that $|\widehat{ \mu^{[0,1)}_{\epsilon}}(\xi)| \leq C_1 (1 + |\xi|)^{-\alpha/2 + \epsilon}$ for all $\xi \in\mathbb{Z}\backslash\lbrace 0 \rbrace$. Let $\phi \in C_c^{\infty}$ be any smooth function supported in $[0,1]$ such that $\phi \mu_{\epsilon}^P$ is not the zero measure. Then there exists a $C_2$ not depending on $\xi$ such that $\widehat{\phi \mu_{\epsilon}^P}(\xi) \leq C_2 (1 + |\xi|)^{-\alpha/2 + \epsilon}$ for all $\xi \in \mathbb{R}$.
\end{mylem}

\begin{proof}
Observe that the Fourier transform of $\mu_{\epsilon}^P$, viewed as a tempered distribution on $\mathbb{R}$, is given by
\[\sum_{s = -\infty}^{\infty} \hat \mu^{[0,1)}(s) \delta_s,\]
where $\delta_s$ is the Dirac mass centered at $s$.

Therefore, we can make sense of $\widehat{\phi \mu_{\epsilon}^P}$ as the convolution of $\widehat{ \mu_{\epsilon}^P}$ and $\hat \phi$. This convolution is equal to
\[ \int \hat \phi(\xi - s) d \widehat{ \mu_{\epsilon}^P} (s)  = \sum_{s = -\infty}^{\infty} \hat \phi(\xi - s) \widehat{ \mu_{\epsilon}^{[0,1)}}(s).\]
We now apply our assumption on $\widehat{ \mu_{\epsilon}^{[0,1)}}(s)$, as well as the Schwartz bound on $\hat \phi$, to conclude
\begin{equation}\label{phimuhat}
|\widehat{\phi \mu_{\epsilon}^{P}}(\xi)| \lesssim \sum_{s = -\infty}^{\infty} (1 + |\xi - s|)^{-100} (1 + |s|)^{-\alpha/2 + \epsilon}.
\end{equation}
We will now write the sum in \eqref{phimuhat} as $S_1 + S_2$, where
\begin{IEEEeqnarray*}{rCl}
S_1 & = & \sum_{|s - \xi| \leq |\xi|/2} (1 + |\xi - s|)^{-100} (1 + |s|)^{- \alpha/2 + \epsilon}.\\
S_2 & = & \sum_{|s - \xi| > |\xi|/2} (1 + |\xi - s|)^{-100} (1 + |s|)^{- \alpha/2 + \epsilon}.\\
\end{IEEEeqnarray*}
We will first estimate $S_1$. Observe that if $|s - \xi| \leq |\xi|/2$, we must have $1 + |s| \geq |\xi|/3$. Therefore, we have
\begin{IEEEeqnarray*}{rCl}
S_1 & \lesssim & (1 + |\xi|/3)^{-\alpha/2 + \epsilon} \sum_{|u| \leq |\xi|/2} (1 + |u|)^{-100} \\
& \lesssim & (1 + |\xi|)^{-\alpha/2 + \epsilon}.
\end{IEEEeqnarray*}
This gives the desired estimate for $S_1$. It remains to estimate $S_2$. In order to estimate $S_2$, observe that the inequality $|s - \xi| \geq |\xi|/2$ implies that $|s - \xi| \geq |s|/10$. Applying this estimate gives
\begin{IEEEeqnarray*}{rCl}
S_2 & \leq & \sum_{|s - \xi| \geq |\xi|/2} (1 + |\xi|/2)^{-99} (1 + |s|/10)^{-1} (1 + |s|)^{-\alpha/2 + \epsilon} \\
& \lesssim & (1 + |\xi|)^{-99} \sum_{|s - \xi| \geq |\xi|/2} (1 + |s|)^{-1 -\alpha/2 + \epsilon} \\
& \leq & (1 + |\xi|)^{-99} \sum_{s = -\infty}^{\infty} (1 + |s|)^{-1 -\alpha/2 + \epsilon} \\
& \lesssim & (1 + |\xi|)^{-99}.
\end{IEEEeqnarray*}
Adding $S_1$ and $S_2$ gives the result.
\end{proof}
We will need one more result that goes in the other direction---a result that allows us to lift compactly supported, bounded, measurable functions $f$ on $\mathbb{R}$ to bounded, measurable functions $f^{[0,1)}$ the torus. We emphasize that the following lemma allows us to control the Fourier coefficients of $f^{[0,1)}$ by knowing $\hat f(s)$ for \textit{integer} values $s$.

Let $f \in L^{\infty}(\mathbb{R})$  be compactly supported. Define $f^P$ by 
\[f^P(x) = \sum_{j \in \mathbb{Z}} f(x + j).\]
The assumptions on $f$ guarantee that $f^P(x)$ converges a.e. to a $1$-periodic function in $L^{\infty}(\mathbb{R})$. This function can naturally be associated to a function $f^{[0,1)}$ on the torus.
\begin{mylem}\label{liftinglemma}
Let $f \in L^{\infty}(\mathbb{R})$ be a compactly supported function, and define $f^{[0,1)}$ as above. Then $\widehat{f^{[0,1)}}(s) = \hat f(s)$ for all integers $s$.
\end{mylem}
\begin{proof}
We have 
\begin{IEEEeqnarray*}{rCl}
\widehat{f^{[0,1)}}(s) & = & \int_0^1 e^{-2 \pi i s x} f^{[0,1)}(x) \, dx \\
		   & = & \int_0^1 e^{-2 \pi i s x} f^{P}(x) \, dx \\
		   & = & \int_0^1 e^{-2 \pi i s x} \sum_{j \in \mathbb{Z}} f(x + j) \, dx \\
		   & = & \sum_{j \in \mathbb{Z}} \int_{j}^{j+1} e^{-2 \pi i s (x - j)} f(x) \, dx \\
		   & = & \sum_{j \in \mathbb{Z}} \int_{j}^{j+1} e^{-2 \pi i s x} f(x) \, dx \\
		   & = & \int_{\mathbb{R}} e^{-2 \pi i s x} f(x) \, dx \\
		   & = & \hat f(s).
\end{IEEEeqnarray*}
\end{proof}
\section{A single-scale estimate}
Lemma \ref{torusmeasure} reduces the proof of Theorem \ref{mainthm2} to finding a measure $\mu$ supported on the torus such that $\hat \mu(s) \leq |s|^{-\alpha/2 + \epsilon}$ for integers $s$. This measure will be constructed as a weak-limit of products of functions, each of which is a sum of smoothed indicator functions of balls of an appropriate scale.

For now, we consider functions supported on the interval $\mathbb{R}$. We will later lift this function to the torus. We define a function $g_{M}$ at scale $M$ which we use to construct our measure supported in the exact order set. The function we consider is supported in the set 
\begin{equation}\label{eq:gMsupp}\left\lbrace \psi_1(q) - c_M\psi_2(q)\leq\left|x-\frac{r-\theta}{q}\right|\leq \psi_1(q);\text{q prime},M\leq q<2M,0\leq r<q \right\rbrace \subset \mathbb{R},\end{equation}
where we take 
\begin{equation}\label{cchoice}
c_M=M^{-\epsilon/100}.
\end{equation} 
Observe that  $0<c_M<1$, with $c_M$ close to $0$ if $M$ is chosen sufficiently large. For the remainder of this section, we will suppress the dependence on $M$ and write $c$ for $c_M$.

For a given prime $q\in [M,2M)$, the interval $I_{q,r}=\{ x : \psi_1(q) - c\psi_2(q)\leq x-\frac{r-\theta}{q}\leq \psi_1(q)\}$ can be expressed as 
\[I_{q,r}=x_{q,r}+c\psi_2(q)[-1/2,1/2],\]
where $x_{q,r}= \frac{r-\theta}{q}+\psi_1(q) - (c/2) \psi_2(q)$. 

Let $\phi$ be a smooth function  with $\supp\phi\subset[-1/2,1/2]$ for which $|\hat{\phi}(\xi)|\leq C\exp \left(-|\xi|^{3/4}\right)$ for large $|\xi|$. Such a function $\phi$ is provided by Ingham \cite{Ingham34}, who in fact constructs a real-valued function whose square satisfies the desired properties.

We define 
\[g_M(x)=\sum_{\substack{M \leq q < 2M \\ \text{q prime}}}\sum_{0 \leq r < q}\phi_{r,q}(x),\] 
where 
\[\phi_{r,q}(x)\coloneqq (c \psi_2(q))^{-1} \phi\left((c \psi_2(q))^{-1} \left(x-x_{q,r}\right)\right).\] 
Observe that the Fourier transform of $\phi_{r,q}$ satisfies
\begin{equation}\label{phift}
\hat{\phi_{r,q}}(s)=e(s x_{q,r})\hat{\phi}(c\psi_2(q)s)
\end{equation}
where $e(x) = e^{2 \pi i x}$. We then set $f_M(x)=g_M(x)/\hat{g_M}(0)$ so that $\hat{f_M}(0)=1$. 
\begin{mylem}\label{fhatestimate}
Let $g_M$, $f_M$ be defined as above. Then 
\begin{IEEEeqnarray}{rCll}
\hat f_M(0) & = & 1 ,\label{fmhat0}\\
\hat f_M(s) & = & 0 \quad & \text{if $1 \leq |s| < M$}, \label{fmhatsmall}\\
|\hat f_M(s)| & \leq & C_{\epsilon} M^{-1 + \epsilon} \quad & \text{if }M \leq |s| \leq M^{\tau (1 + \epsilon/2)}, \label{fmhatmedium}\\
|\hat f_M(s)| & \leq & \exp \left( -|\frac{s}{M^{\tau}}|^{1/2} \right) \quad & \text{if }|s| \geq M^{\tau (1 + \epsilon/2)} \label{fmhatlarge}.
\end{IEEEeqnarray}
\end{mylem}
\begin{proof}
Equation \eqref{fmhat0} follows directly from our normalization of $f_M$.

For the proof of \eqref{fmhatsmall}, \eqref{fmhatmedium}, and \eqref{fmhatlarge}, we will begin by computing $\hat g_M(s)$ explicitly. By equation \eqref{phift}, we have 
\begin{equation}\label{gmhat}
\hat g_M(s) = \sum_{\substack{M \leq q < 2M \\ q \text{ prime}}} \sum_{0 \leq r < q} e(sx_{q,r}) \hat{\phi}(c\psi_2(q)s).
\end{equation}
By plugging in the value for $x_{q,r}$, we obtain
\begin{IEEEeqnarray*}{rCl}
\hat g_M(s) & = & \sum_{\substack{M \leq q < 2M \\ q \text{ prime}}} \sum_{0 \leq r < q} e \left(s \left( \frac{r-\theta}{q}+ \psi_1(q) - (c/2) \psi_2(q) \right) \right) \hat{\phi}(\psi_2(q)s) \\
& = & \sum_{\substack{M \leq q < 2M \\ q \text{ prime}}} e \left(s \left( \frac{- \theta}{q} + \psi_1(q) - (c/2) \psi_2(q) \right) \right) \hat{\phi}(c \psi_2(q)s)\sum_{0 \leq r < q} e \left(\frac{rs}{q}\right).
\end{IEEEeqnarray*}
The inner sum $\sum_{0 \leq r < q} e \left(\frac{rs}{q} \right) $ is a geometric series that evaluates to $0$ unless $q | s$, in which case it evaluates to $q$. Therefore, we have
\begin{IEEEeqnarray*}{rCl}
\hat g_M(s) & = &  \sum_{\substack{M \leq q < 2M \\ q \text{ prime} \\ q | s}} e \left(s \left( \frac{- \theta}{q} + \psi_1(q) - (c/2) \psi_2(q) \right) \right) \hat{\phi}(c \psi_2(q) s) \cdot q.
\end{IEEEeqnarray*}
For $0 < |s| < M$ this sum is empty, establishing that $\hat g_M(s) = 0$ and giving \eqref{fmhatsmall}. If $s = 0$, the sum in $q$ consists of all $M \leq q < 2M$, giving
\begin{equation}\label{gmhat0}
\hat g_M(0) =  \hat \phi(0) \sum_{\substack{M \leq q < 2M \\ q \text{ prime}}} q \sim \frac{M^2}{\log M}
\end{equation}
by the prime number theorem. 

For other values of $s$, we use the triangle inequality to give the estimate
\begin{equation}\label{gmhatmedlarge}
|\hat g_M(s)| \lesssim M \sum_{\substack{M \leq q < 2M \\ q \text{ prime} \\ q | s}} |\hat{\phi}(c\psi_2(q)s)|.
\end{equation}
Combining \eqref{gmhat0} and \eqref{gmhatmedlarge} gives the estimate
\begin{equation}\label{fmhatmedlarge}
|\hat f_M(s)| \lesssim M^{-1} \log M \sum_{\substack{M \leq q < 2M \\ q \text{ prime} \\ q | s}} |\hat{\phi}(c\psi_2(q)s)|.
\end{equation}

We now consider the regime where $M \leq |s| < M^{\tau_2 (1 + \epsilon/2)}$. Observe that the number of terms $q$ in the sum on the right hand side of \eqref{fmhatmedlarge} is bounded above by $\frac{\log |s|}{\log M}$. Using the bound $|\hat \phi (c \psi_2(q) s)| \lesssim 1$, we find
\[|\hat f_M(s)| \lesssim M^{-1} \log |s| \lesssim M^{-1} \log M \leq C_{\epsilon} M^{-1 + \epsilon},\]
establishing \eqref{fmhatmedium}.

For $|s |> M^{\tau_2(1 + \epsilon/2)}$, we take advantage of the choice of $\phi$. First, observe that, because $\psi_2$ is decreasing, we can bound $c \psi_2(q)s$ from below by $c\psi_2(2M)s$. By \eqref{cchoice}, this is $M^{-\epsilon/100} \psi_2(2 M)s.$ Provided that $M$ is sufficiently large depending on $\epsilon$, we can use \eqref{orderdef} to conclude that 
\[|c \psi_2(q) s| \geq 2^{-\tau_2 - \epsilon/100} M^{-\tau_2 -\epsilon/50} |s|.\]
Therefore, applying our assumption on $\phi$, we have that
\begin{equation}\label{singleterm}
\hat \phi(c \psi_2(q) s)  \leq C \exp \left(-2^{-3/4(\tau_2 + \epsilon/100)}M^{-3/4(\tau_2 + \epsilon/50)}|s|^{3/4} \right).
\end{equation}
Now, we estimate the series \eqref{fmhatmedlarge} using the bound \eqref{singleterm}. There are no more than $M$ terms in this sum, so, for an appropriate constant $C$,
\[|\hat f_M(s)| \leq C \log M \exp \left(-2^{-3/4(\tau_2 + \epsilon/100)}M^{-3/4(\tau_2 + \epsilon/50)}|s|^{3/4} \right).\]
Since we are in the regime $|s| > M^{\tau_2(1 + \epsilon/2)},$
\[|\hat f_M(s)| \leq C \log M \exp \left(- \left(|s|^{1/2} M^{-\tau_2/2} \right) \left(2^{- 3/4 (\tau_2 + \epsilon/100)} M^{\epsilon \tau_2/8 - 3 \epsilon/200} \right) \right). \]
Because the exponent $\epsilon \tau_2/8 - 3 \epsilon/200$ is positive, it follows that for $M$ sufficiently large, we have the bound
\[|\hat f_M(s)| \leq \exp \left(-|s|^{1/2} M^{-\tau_2/2} \right), \]
establishing the desired bound \eqref{fmhatlarge}.
\end{proof}
\section{A Convolution Stability Lemma}
In this section, we establish a convolution stability lemma. This lemma will later be used in Section \ref{sec:construction}, in combination with Lemma \ref{fhatestimate}, applied at different scales as part of an induction argument, to complete the proof of Theorem \ref{mainthm2}.

To this end, we will consider a sequence $\{M_j\}_{j=1}^{\infty}$ of positive numbers whose growth rate is dictated by Lemma \ref{DiophantineLemmaGeneral}.

This convolution stability lemma provides an estimate for $F * G$, where functions $F$ and $G$ satisfy certain bounds following Lemma \ref{fhatestimate}. In practice, the function $G$ will be $\hat f_{M_1} * \hat f_{M_2} * \cdots * \hat f_{M_{j}}$ for some appropriate $j$, and $F$ will be taken to be $\hat f_{M_{j+1}}$. 

In this section, we will assume that $\tau_1$ and $\tau_2$ satisfy the condition $\beta > \tau_2$. Recall we defined $\beta_\epsilon=\gamma^{-2}(\tau_1- 1)^2 - \epsilon=\beta-\epsilon$ in the equation \eqref{eq:betaEps}. We consider only those $\epsilon$ small enough so that $\beta_\epsilon>\tau_2$.

\begin{mylem}[Convolution Stability Lemma]\label{convstab}
Let $M_{j+1} = M_j^{\beta_{\epsilon}}$, where $\beta_{\epsilon} > \tau_2 \geq 2$. Let $\epsilon > 0$ be sufficiently small that the quantity
\begin{equation}\label{deltachoice}
\delta : = \frac{\beta_{\epsilon} (1 - \epsilon) - \tau_2 (1 + \epsilon)^2}{\tau_2(\beta_{\epsilon} - 1)(1 + \epsilon)},
\end{equation}
satisfies $\delta - \epsilon > 0$, and let $F, G : \mathbb{Z} \to \mathbb{C}$ be functions satisfying the following estimates:
\begin{IEEEeqnarray}{rCll}
F(0)  & = & 1  \label{eq:Fat0},\\
F(s) & = & 0 \quad & \text{if $1 \leq |s| < M_{j+1}$} \label{eq:FsSmll},\\
|F(s)| & \leq &  C_{\epsilon} M_{j+1}^{-1 + \epsilon} & \text{if $M_{j+1} \leq |s| \leq M_{j+1}^{\tau_2 (1 + \epsilon/2)}$}  \label{eq:FsDecayIntrmdt},\\
|F(s)| & \leq & \exp \left( -\left|\frac{s}{M_{j+1}^{\tau_2}}\right|^{1/2} \right) \quad & \text{if }|s| \geq M_{j+1}^{\tau_2 (1 + \epsilon/2)} \label{eq:FsDecayLrg}.
\end{IEEEeqnarray}
and
\begin{IEEEeqnarray}{rCll}
G(0)  & \leq & 2  \label{eq:Gat0},\\
|G(s)| & \leq & 2 |s|^{-\delta + \epsilon} & \text{if $|s| \leq M_j^{\tau_2 (1 + \epsilon)}$}  \label{eq:GsDecayIntrmdt},\\
|G(s)| & \leq &\exp \left( - \frac{1}{2} \left| \frac{s}{M_j^{\tau_2}}\right|^{1/2} \right) & \text{if $|s| \geq M_j^{\tau_2(1 + \epsilon)}$} \label{eq:GsDecayLrg}.
\end{IEEEeqnarray}
Then, provided that $M_1$ is sufficiently large depending on $\epsilon$, we have the following three conclusions:

\begin{enumerate}[(a)]
\item \[|F * G(s) - G(s)| \leq M_{j+1}^{-\delta} \quad \text{if $|s| \leq M_j^{\tau_2(1 + \epsilon)}$}\] \label{smalls} 
\item \[|F * G(s)| \leq |s|^{- \delta + \epsilon} \quad  \text{if $M_j^{\tau_2 (1 + \epsilon)} < |s| \leq M_{j+1}^{\tau_2(1 + \epsilon)}$} \] \label{mediums}
\item \[|F * G(s)| \leq \exp \left(- \frac{1}{2} \left| \frac{s}{M_{j+1}^{\tau_2}} \right|^{1/2} \right) \quad  \text{if $|s|> M_{j+1}^{\tau_2(1 + \epsilon)}$}.\] \label{larges}
\end{enumerate}
\end{mylem}

For reference, we here give a sketch of $F$ and $G$ corresponding to some index $j$. The usage of $\approx$ in this figure indicates an $\epsilon$-loss in the exponent on $M_j$ or $M_{j+1}$.

\begin{figure}[h]
\includegraphics[width=\textwidth]{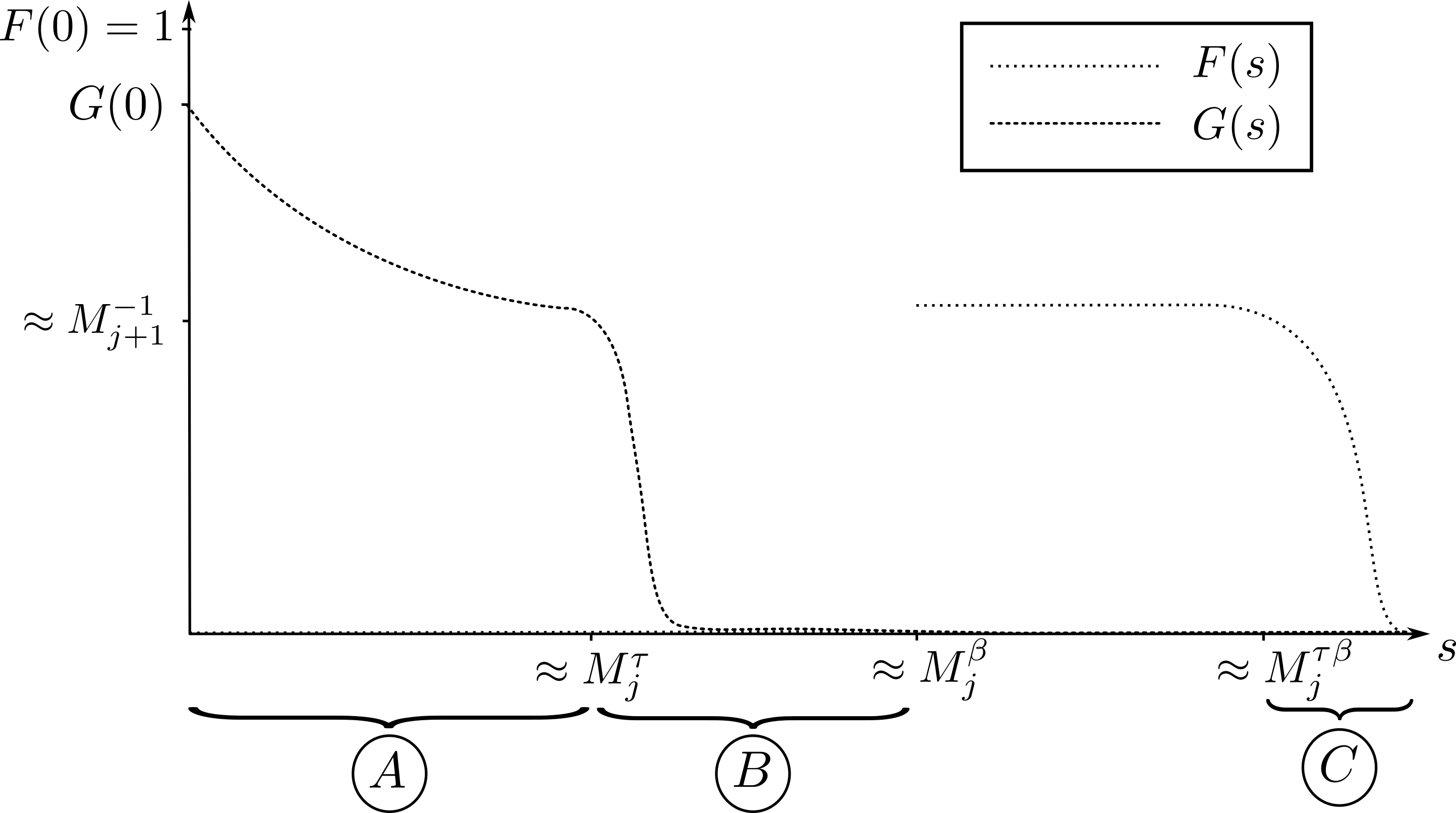}
\centering
\caption{Sketch of $F$ and $G$}
\end{figure}

Note that in region $A$, we have that $|G(s)|\lesssim|s|^{-\delta}$. In region $B$, $G(s)$ decays rapidly. In regions $A$ and $B$ (not including $0$), $F$ vanishes. In region $C$, $F$ decays rapidly.
\begin{myrmk}
Lemma \ref{convstab} is called a \textbf{convolution stability lemma} because of the bound (\ref{smalls}), which shows that for small values of $s$, the convolution $F * G$ will be very close to $G$.

The bound (\ref{mediums}) will dictate the Fourier decay of the infinite product measure supported on our set. Note that the bound $|s|^{- \delta + \epsilon}$ is significantly worse than the bound $M_{j+1}^{-1 + \epsilon}$ available for $F$ in this region- this is the reason for the loss in Fourier dimension compared to the set of well-approximable numbers from Kaufman's argument.

The bound (\ref{larges}) will allow for the convolution stability lemma to be applied inductively. Although this bound gets slightly worse at each stage of the induction, it will always be good enough to match the conditions required for $G$ at the next stage of the induction.

\end{myrmk}
\begin{proof}
We will first prove (\ref{smalls}). To this end, we assume $|s| \leq M_j^{\tau_2(1 + \epsilon)}$. We write
\[F * G(s) = \sum_{t \in \mathbb{Z}} F(t) G(s - t).\]
The main contribution to this sum will come from the $t = 0$ term, which is precisely $G(s)$. Additionally, there is no contribution for $1 \leq |t| < M_{j+1}$ because $F(t) = 0$ there. Thus, we see that 
 \begin{equation}\label{eq:convStabSmallsBound}\left|F*G(s)-G(s)\right|\leq  \left|\sum_{M_{j+1}\leq |t|\leq M_{j+1}^{\tau_2(1 + \epsilon)}} F(t) G(s - t) \right|+  \left|\sum_{M_{j+1}^{\tau_2(1 + \epsilon)}<|t |} F(t) G(s - t)\right|.\end{equation}

We now estimate the size of the first term on the right hand side of \eqref{eq:convStabSmallsBound} and consider the corresponding $M_{j+1}\leq |t|\leq M_{j+1}^{\tau_2(1 + \epsilon)}$. We have $|F(t)| \leq C_{\epsilon} M_{j+1}^{-1 + \epsilon}$. Furthermore, because $|s| \leq M_j^{\tau_2 (1 + \epsilon)}$,  we have $|s - t| \geq \frac{M_{j+1}}{2}$, provided $M_1$ is chosen large enough.
 Thus, we see 
\begin{IEEEeqnarray*}{rCl}
|G(s- t)| & \leq & \exp \left(- \frac{1}{2} \left(\frac{M_{j+1}}{2 M_j^{\tau_2}}\right)^{1/2} \right)\\
& \leq & \exp \left( - \frac{1}{2} \frac{M_j^{\frac{\beta_{\epsilon} - \tau_2}{2}}}{\sqrt{2}} \right).
\end{IEEEeqnarray*}
Recall that $\beta_\epsilon-\tau_2>0$. Combining the bounds on $F$ and $G$ and counting the number of terms in the sum, the first sum of \eqref{eq:convStabSmallsBound} is bounded by
\[C_{\epsilon} \exp \left( - \frac{1}{2} \frac{M_j^{\frac{\beta_{\epsilon} - \tau_2}{2}}}{\sqrt{2}} \right) M_{j+1}^{-1 + \epsilon + \tau_2(1 +\epsilon)}= C_{\epsilon} \exp \left( - \frac{1}{2} \frac{M_{j+1}^{\frac{\beta_{\epsilon} - \tau_2}{2\beta_{\epsilon}}}}{\sqrt{2}} \right) M_{j+1}^{-1 + \epsilon + \tau_2(1 +\epsilon)}\leq \frac{1}{2} M_{j+1}^{-\delta},\]
provided $M_1$ is chosen large enough.

The final step in proving (\ref{smalls}) is to bound the second term on the right hand side of \eqref{eq:convStabSmallsBound}. For $|t|>M_{j+1}^{\tau_2(1 + \epsilon)}$, we still have $|s - t| \geq M_{j}^{\tau_2(1 + \epsilon)}$, so we can apply the tail estimates for both $F$ and $G$. Thus,
\[F(t) \leq \exp \left( - \left|\frac{t}{M_{j+1}^{\tau_2}}\right|^{1/2} \right).\] For simplicity, we observe that $|s - t| \geq |t|/2$ and of course $|G(s - t)| \leq 1$.  So we estimate 
\[ \left|\sum_{M_{j+1}^{\tau_2(1 + \epsilon)}<|t| } F(t) G(s - t)\right|\leq\sum_{|t| \geq M_{j+1}^{\tau_2(1 + \epsilon)}} \exp \left(-\left|\frac{t}{M_{j+1}^{\tau_2}}\right|^{1/2} \right)\leq \frac{1}{2} M_{j+1}^{-\delta},\]
by comparing to the corresponding integral, if $M_1$ is sufficiently large. Summing the two terms completes the proof of (\ref{smalls}).

We now prove (\ref{mediums}). We consider those $s$ such that $M_{j}^{\tau_2(1 + \epsilon)} < |s| \leq M_{j+1}^{\tau_2(1 + \epsilon)}$. We bound $F*G$,
\[\left|\sum_{t \in \mathbb{Z}} F(s - t) G(t)\right|\]
\begin{equation}\label{eq:convStabMedsBound}\leq \left|\sum_{|t|\leq M_j^{\tau_2(1 + \epsilon)} } F(s - t) G(t)\right|+ \left|\sum_{M_{j}^{\tau_2(1 + \epsilon)} \leq |t| \leq 2 M_{j+1}^{\tau_2(1 + \epsilon)} } F(s - t) G(t)\right|+\left|\sum_{2M_{j}^{\tau_2(1 + \epsilon)} \leq |t| } F(s - t) G(t)\right|.\end{equation}
The main contribution to this bound will be the first term. For such values of $|t|\leq M_j^{\tau_2(1 + \epsilon)}$, we have the estimate $|G(t)| \leq 2 |t|^{-\delta + \epsilon}$ for $t \neq 0$ and $|G(t)| \leq 2$ at $t = 0$. 

Continuing our analysis of the first term of \eqref{eq:convStabMedsBound}, where $|t|\leq M_j^{\tau_2(1 + \epsilon)}$, we have that $|F(s - t)| \leq C_{\epsilon} M_{j+1}^{-1 + \epsilon}$. As a result, for an appropriate constant $K$,
\[\left|\sum_{|t|\leq M_j^{\tau_2(1 + \epsilon)} } F(s - t) G(t)\right|\leq M_{j+1} ^{- 1 + \epsilon} \left(2C_{\epsilon} + 4 C_{\epsilon} \sum_{t=1}^{M_j^{\tau_2(1 + \epsilon)}} |t|^{-\delta + \epsilon} \right).\]
\[\leq KM_{j+1}^{-1 + \epsilon + (- \delta + \epsilon + 1) \frac{\tau_2}{\beta_{\epsilon}}(1 + \epsilon)} \leq K |s|^{\frac{1}{\tau_2 (1 + \epsilon)} (-1 + \epsilon + (- \delta + \epsilon + 1) \frac{\tau_2}{\beta_{\epsilon}} (1 + \epsilon))}\leq\frac{1}{3} |s|^{-\delta + \epsilon},\] 
since our choice of $\delta$ guarantees that the penultimate exponent on $|s|$ above is $- \delta$ and we can absorb the constant $K$ in to the $|s|^{\epsilon}$ term, which is possible if $M_1$ is sufficiently large.

We now consider the second term in the bound \eqref{eq:convStabMedsBound}. The relevant $t$ are those with $M_{j}^{\tau_2(1 + \epsilon)} \leq |t| \leq 2 M_{j+1}^{\tau_2(1 + \epsilon)}$, including the case in which $t = s$. For such $t$, we have a bound of $1$ on $|F(s - t)|$ and a bound of $\exp \left( - \frac{1}{2} \left|\frac{t}{M_j^{\tau_2}} \right|^{1/2} \right)$ for $|G(t)|$. As there are at most $4M_{j+1}^{\tau_2(1 + \epsilon)}$ such values of $t$, 
\[  \left|\sum_{M_{j}^{\tau_2(1 + \epsilon)} \leq |t| \leq 2 M_{j+1}^{\tau_2(1 + \epsilon)} } F(s - t) G(t)\right|\leq 4 M_{j+1}^{\tau_2 (1 + \epsilon)} \exp \left( - \frac{1}{2} \left| \frac{M_j^{\tau_2(1 + \epsilon)}}{M_j^{\tau_2}} \right|^{1/2} \right)\leq\frac{1}{3} |s|^{- \delta + \epsilon},\] if $M_1$ is taken large enough depending on the choice of $\epsilon$.

Finally, we consider the final term in \eqref{eq:convStabMedsBound}. Here $|t| \geq 2 M_{j+1}^{\tau_2 (1 + \epsilon)}$, and we are well within the region on which the tail bounds can be applied for both $F$ and $G$. We will not need both tail bounds, however; we will simply bound $|F(s - t)|$ by $1$ on this region, and use the tail bound $\exp \left( - \frac{1}{2} \left|\frac{t}{M_j^{\tau_2}} \right|^{1/2} \right)$ for $G(t)$. Then, by comparing to the integral and ensuring that $M_1$ is sufficiently large, we find \[\left|\sum_{2M_{j}^{\tau_2(1 + \epsilon)} \leq |t| } F(s - t) G(t)\right|\leq \frac{1}{3} |s|^{-\delta + \epsilon}.\]

The established bounds on each of the terms in \eqref{eq:convStabMedsBound} combine to show $|F * G(s)| \leq |s|^{-\delta + \epsilon}$, completing the proof of (\ref{mediums}). 

It remains to prove (\ref{larges}).
Let $s$ be such that $|s| > M_{j+1}^{\tau_2(1 + \epsilon)}$. Writing the convolution as for (\ref{mediums}), we bound $F * G(s)$ as follows:
\[\left| \sum_{t \in \mathbb{Z}} F(s - t) G(t) \right|\]
\begin{equation}\label{eq:convStabLargesBound} \leq \left|\sum_{\substack{|t|\leq 2|s| \\ |s-t|\geq \frac{|s|}{2}}}F(s - t) G(t)\right|+\left|\sum_{\substack{|t|\leq 2|s| \\ |s-t|< \frac{|s|}{2} }}F(s - t) G(t)\right|+\left|\sum_{|t|> 2|s| }F(s - t) G(t)\right|.\end{equation}

The thrust of the proof is that, because $|s|$ is so large, we are always in a situation for which the tail bounds on either $F$ or $G$ will apply. 

We consider the first term of \eqref{eq:convStabLargesBound}, where $|t| \leq 2 |s|$ and $|s - t| \geq \frac{|s|}{2}$. For such $t$, including $t = 0$, we have the bound $|G(t)| \leq 2$.  On the other hand, for such $t$, we certainly have $|s - t| \geq \frac{1}{2} M_{j+1}^{\tau_2(1 + \epsilon)} \geq M_{j+1}^{\tau_2(1 + \epsilon/2)}$, provided $M_1$ is taken sufficiently large. Thus, we have
\[|F(s - t)| \leq \exp \left( - \left|\frac{s}{2 M_{j+1}^{\tau_2}} \right|^{1/2} \right).\]
Summing over $|t|\leq 2|s|$, we see
\[\left|\sum_{\substack{ |t|\leq 2|s|  \\|s-t|\geq \frac{|s|}{2}}}F(s - t) G(t)\right|
\leq 4 |s| \exp \left(- \left| \frac{s}{2 M_{j+1}^{\tau_2}} \right|^{1/2} \right)
\leq\frac{1}{3} \exp \left( - \frac{1}{2} \left| \frac{s}{M_{j+1}^{\tau_2}} \right|^{1/2} \right),\]
provided $M_1$ is large enough depending on $\epsilon$.

Next, we will bound the second term of \eqref{eq:convStabLargesBound}. The relevant $t$ are such that $|s - t| < \frac{|s|}{2}$ (including the $t = s$ term). Note that for such $t$, we certainly have $|t| \geq \frac{|s|}{2} \geq M_j^{\tau_2(1 + \epsilon)}$, if $M_1$ is large enough. We also have $|F(s - t)| \leq 1$, and $|G(t)| \leq \exp \left( - \frac{1}{2} \left| \frac{s}{2M_j^{\tau_2}} \right|^{1/2} \right)$. Observe that the total number of values of $t$ summed is at most $4|s|$. If $M_1$ is sufficiently large, keeping in mind that $|s| \geq M_{j+1}^{\tau_2(1 + \epsilon)}$, we observe that $4|s|$ is much less than $\frac{1}{3} \exp( \frac{1}{2} |s|^{1/2} ((2M_j)^{-\tau_2/2} - M_{j+1}^{-\tau_2/2}))$. Thus,
\[\left|\sum_{\substack{|t|\leq 2|s| \\ |s-t|< \frac{|s|}{2} }}F(s - t) G(t)\right|\leq 4|s|\exp \left( - \frac{1}{2} \left| \frac{s}{2M_j^{\tau_2}} \right|^{1/2} \right)\leq\frac{1}{3} \exp \left(- \frac{1}{2} \left| \frac{s}{M_{j+1}^{\tau_2}} \right|^{1/2} \right).\]

It remains to bound the third term in \eqref{eq:convStabLargesBound}.  Here, $|t| \geq 2|s|$, and we have $|s - t| \geq \frac{|t|}{2}$. We will use an estimate of $1$ for $|G(t)|$ and a bound of $\exp \left( - \left| \frac{t}{2 M_{j+1}^{\tau_2}} \right|^{1/2} \right)$ for $|F(s-t)|$, as we can apply the tail bound on $F$. Hence,
\[\left|\sum_{|t|> 2|s| }F(s - t) G(t)\right|\leq\sum_{|t| \geq 2 |s|} \exp \left(- \left| \frac{t}{2 M_{j+1}^{\tau_2}} \right|^{1/2} \right)\leq\frac{1}{3} \exp \left( - \frac{1}{2} \left| \frac{s}{M_{j+1}^{\tau_2}} \right|^{1/2} \right),\]
by comparison with the corresponding integral.

We arrive at the desired bound, (\ref{larges}), by summing the three terms.
\end{proof}

\section{Construction of the Measure}\label{sec:construction}
In this section, we complete the proof of Theorem \ref{mainthm2}. Recall that \[\delta = \frac{\beta_{\epsilon} (1 - \epsilon) - \tau_2 (1 + \epsilon)^2}{\tau_2(\beta_{\epsilon} - 1)(1 + \epsilon)};\] this was defined at \eqref{deltachoice}. The stated lower bound on the Fourier dimension, $\alpha= \frac{2(\beta - \tau_2)}{\tau_2 (\beta - 1)}$, was given at \eqref{fourierdimestimate2}. To prove Theorem \ref{mainthm2}, it suffices to construct a measure $\mu_{\epsilon}$ on the torus $[0,1)$, with the decay $|\hat \mu_{\epsilon}(s)|\lesssim |s|^{-\delta+\epsilon}$.

Let $\epsilon > 0$, and let $M_1$ be a number so large that Lemma \ref{convstab} applies (with $M_{j+1} = M_{j}^{\beta_{\epsilon}}$ for all $j > 1$), and sufficiently large that $\sum_{j=1}^{\infty} M_j^{-\delta}  = \sum_{j=1}^{\infty} M_1^{-j \beta_{\epsilon} \delta} < \frac{1}{100}$. For each $j$, define $f_{M_j}$ as in Lemma \ref{fhatestimate}, and $f_{M_j}^{[0,1)}$ as in Lemma \ref{liftinglemma}. We define the function $\mu_{\epsilon}^{(k)} = \prod_{j=1}^{k} f_{M_j}^{[0,1)}$. We claim that the measures $\mu_{\epsilon}^{(k)}$ have a subsequence with a weak limit $\mu_{\epsilon}$ with the desired properties. In the proof, we conflate the function $\mu_{\epsilon}^{(k)}$ with the absolutely continuous measure whose Radon-Nikodym derivative is $\mu_{\epsilon}^{(k)}$.

The proof of this will require us to estimate $\hat \mu_{\epsilon}^{(k)}(s)$ for integer values $s$. We will obtain the following estimate by applying Lemma \ref{convstab} inductively. 
\begin{mylem}\label{inductiveestimate}
Let $M_0 = 0$ for convenience. We have the following estimates on $\hat \mu_{\epsilon}^{(k)}(s)$ for any integers $k \geq 1$ and $s \in \mathbb{Z}$:
\begin{equation}\label{indest0}
1 - \sum_{j=1}^{k} M_j^{-\delta} < \hat \mu_{\epsilon}^{(k)}(0) < 1 + \sum_{j=1}^k M_j^{-\delta},  \\
\end{equation}
\begin{equation}\label{indestmed}
|\hat \mu_{\epsilon}^{(k)}(s)| < |s|^{-\delta + \epsilon} + \sum_{j=J}^{k} M_{j}^{-\delta} < 2 |s|^{-\delta + \epsilon} \quad \text{if $M_{J - 1}^{\tau_2(1 + \epsilon)} < |s| \leq M_J^{\tau_2(1 + \epsilon)}$ for $1 \leq J \leq k$},
\end{equation}
\begin{equation}\label{indestlarge}
|\hat \mu_{\epsilon}^{(k)}(s)| \leq \exp \left(- \frac{1}{2} \left| \frac{s}{M^{\tau_2}_{k}} \right|^{1/2} \right) \quad \text{if $|s|> M_{k}^{\tau_2(1 + \epsilon)}$}.
\end{equation}
\end{mylem}
\begin{proof}
We prove this lemma by induction. The base case of this lemma is implied by Lemma \ref{fhatestimate} applied to $f_{M_1}^{[0,1)}$. So we need only show the inductive step.

Suppose, for some $k > 1$, we have that $\mu_{\epsilon}^{(k-1)}$ satisfies the estimates in Lemma \ref{inductiveestimate}. We must show that $\mu_{\epsilon}^{(k)}$ also satisfies these estimates. Our tool for this is Lemma \ref{convstab}. Observe that, by definition, we have that 
\[\mu_{\epsilon}^{(k)} = \mu_{\epsilon}^{(k-1)} f_{M_k}^{[0,1)}.\] 
Therefore, by the convolution rule for the Fourier transform, we have
\[\hat \mu_{\epsilon}^{(k)} = \hat \mu_{\epsilon}^{(k-1)} * \widehat{f_{M_k}^{[0,1)}}.\]
The estimates \eqref{indest0}, \eqref{indestmed}, and \eqref{indestlarge} imply that $\hat \mu_{\epsilon}^{(k-1)}$ is able to serve as the function $G$ in Lemma \ref{convstab}. Note that \eqref{indest0} for $\mu_{\epsilon}^{(k)}$ follows immediately by combining (\ref{smalls}) of Lemma \ref{convstab} and \eqref{indest0} for $\mu_{\epsilon}^{(k-1)}$. Similarly, for $1 \leq J \leq k-1$,  we have that \eqref{indestmed} holds for $\mu_{\epsilon}^{(k)}$ by combining (\ref{smalls}) of Lemma \ref{convstab} and \eqref{indestmed} for $\mu_{\epsilon}^{(k-1)}$. The $J = k$ case of the estimate \eqref{indestmed} for $\mu_{\epsilon}^{(k)}$ is an immediate consequence of (\ref{mediums}) of Lemma \ref{convstab}. Finally, the estimate \eqref{indestlarge} is given by (\ref{larges}) of Lemma \ref{convstab}.
\end{proof}
We are now in a position to define our measure $\mu_{\epsilon}$. It is clear from the Banach-Alaoglu theorem that some subsequence of the $\mu_{\epsilon}^{(k)}$ converges weakly to some measure $\mu_{\epsilon}$. The estimate \eqref{indest0} shows that the weak-limit of this subsequence is a nonzero finite measure, and it is clear from the fact that $\mu_{\epsilon}^{(k)} \geq 0$ for all $k$ that $\mu_{\epsilon} \geq 0$. The estimate \eqref{indestmed} implies that $|\hat \mu_{\epsilon}(s)| \leq 2 |s|^{-\delta + \epsilon}$ for all $s$. Therefore, in order to establish the Fourier dimension bound, the only statement it remains to prove about $\mu_{\epsilon}$ is that its support is contained in $\Tight^{[0,1)}(\psi_1, \psi_2, \theta)$.
\begin{mylem}\label{muepsilonsupport}
The measure $\mu_{\epsilon}$ is supported on $\Tight^{[0,1)}(\psi_1, \psi_2, \theta)$.
\end{mylem}
\begin{proof}
The measures $\mu_{\epsilon}^{(k)}$ have nested, decreasing support, so we must have 
\[\supp \mu_{\epsilon} \subset \bigcap_k \supp \mu_{\epsilon}^{(k)} = \bigcap_{k} \supp f_{M_k}^{[0,1)}.\]
So it is sufficient to prove that
\[\bigcap_k \supp f_{M_k}^{[0,1)} \subset \Tight^{[0,1)}(\psi_1, \psi_2, \theta).\]
By construction, $\supp f_{M_k}$ is contained in the set given at \eqref{eq:gMsupp}, with $c_{M}=c_{M_k}=M_k^{-\epsilon/100}$. Let $c < 1$. There exists $K = K(\epsilon, c)$ such that $M_K^{- \epsilon/100} < c$. 

Suppose $x \in \bigcap_k \supp f_{M_k}^{[0,1)}.$ Let $x^*$ in $\mathbb{R}$ be the element of the interval $[0,1)$ that is congruent modulo $1$ to $x$. Then for any $k$, there exists an integer $z_k$ such that $x^* + z_k \in \supp f_{M_k}$. This means that, for any $k > K$, there exists a pair $(r_k, q_k)$ with $0 \leq r_k < q_k$ and $M_k \leq q_k < 2M_k$ and $q_k$ prime such that $\psi_1(q) - c \psi_2(q) < \left|x^* - \frac{r_k - z_k q_k - \theta}{q_k} \right| < \psi_1(q)$. Letting $r_k' = r_k - z_k q_k$, this gives, for every $k > K$, a pair $(r_k', q_k)$ such that $\psi_1(q_k) - c \psi_2(q_k) < \left|x^* - \frac{r_k'- \theta}{q_k} \right| < \psi_1(q_k)$. Since $q_{k+1} \lesssim q_{k}^{\beta_{\epsilon}}$ for every $k$, Lemma \ref{DiophantineLemmaGeneral} shows that, provided $k > K'$ for an appropriate value $K'(\epsilon, c)$, there no pair $(r,q)$ such that $q_k < q < q_{k+1}$ and $\left|x^* -  \frac{r - \theta}{q} \right| \leq \psi_1(q) - c \psi_2(q)$. Because this works for all $k > K'' := \max(K, K')$, this shows that $\left|x^* -  \frac{r- \theta}{q} \right|> \psi_1(q) - c \psi_2(q)$ for all pairs $(r,q)$ with $q > q_{K''}$, establishing the result.
\end{proof}
We have shown the support of $\mu_{\epsilon}^{(k)}$ is contained in $\Tight^{[0,1)}(\psi_1, \psi_2, \theta)$.  Applying Lemma \ref{torusmeasure} gives the desired measure on $\Tight(\psi_1, \psi_2, \theta)$.
\section*{Acknowledgements}

The authors would like to thank Sanju Velani and Evgeniy Zorin, without whom this project would not have been possible.

This material is based upon work supported by the National Science Foundation under Award No. 1803086.

The second author was supported by The Maxwell Institute Graduate School in Analysis and its Applications, a Centre for Doctoral Training funded by the UK Engineering and Physical Sciences Research Council (Grant EP/L016508/01), the Scottish Funding Council, Heriot-Watt University and the University of Edinburgh.
\bibliographystyle{myplain}
\bibliography{Fourier_Dimension_Exact_Order}
\end{document}